\documentclass[11pt,a4paper]{article}
\pdfoutput=1
\usepackage[T1]{fontenc}
\usepackage[utf8]{inputenc}
\usepackage{lmodern}
\usepackage[english]{babel}
\usepackage[margin=2.4cm]{geometry}
\usepackage{microtype}
\usepackage{amsmath,amssymb,amsthm}
\usepackage{booktabs,array,tabularx}
\usepackage{enumitem}
\usepackage[authoryear,round]{natbib}
\usepackage{tikz}
\usetikzlibrary{calc,arrows.meta}
\usepackage{xcolor}
\usepackage[hidelinks]{hyperref}
\usepackage{url}

\hypersetup{
  pdftitle={Two teams that cancel: a didactic analysis of a teaching progression for signed numbers from Grade 5 to Grade 7},
  pdfauthor={Vincent Tariel-Adam},
  pdfkeywords={signed numbers, integers, negative numbers, zero pair, counter model, number line, sign rule, didactics of mathematics, curriculum design}
}

\newtheorem{principle}{Design principle}
\theoremstyle{definition}
\newtheorem{hypothesis}{Hypothesis}

\usepackage{ifthen}
\colorlet{colordef}{red!60}
\definecolor{colorprop}{rgb}{0.0,0.4,0.74}
\newcommand\tribuNegative[2]{%
\begin{scope}[shift={(#1,#2)}]
\begin{scope}[scale=0.75]
\draw[fill=colorprop, thick] (0,0)circle(0.4);
\draw[ thick] (-0.15,0)--(0.15,0);
\end{scope}
\end{scope}
}
\newcommand\tribuPositive[2]{%

\begin{scope}[shift={(#1,#2)}]
\begin{scope}[scale=0.75]
\draw[fill=colordef, thick] (0,0)circle(0.4);
\draw[ thick] (0,-0.15)--(0,0.15);
\draw[ thick] (-0.15,0)--(0.15,0);
\end{scope}
\end{scope}
}

\newcommand\explosionForce[2]{
\begin{scope}[shift={(#1,#2)}]
  \shade[inner color=colorprop, outer color=colordef] (0,0) circle (0.3);
 \end{scope}
}
\newcommand\additionRelativeNumber[2]{
\pgfmathsetmacro{\numberFirst}{#1}
\pgfmathsetmacro{\numberFirstAbs}{int(abs(\numberFirst)}

\pgfmathsetmacro{\numberSecond}{#2}
\pgfmathsetmacro{\numberSecondAbs}{int(abs(\numberSecond)}

\pgfmathparse{ifthenelse(\numberFirst*\numberSecond>0,1,0)} 
\ifthenelse{\pgfmathresult=1}
{

\foreach \x in {1,...,\numberFirstAbs}  {
\pgfmathparse{ifthenelse(\numberFirst>0,1,0)} \ifthenelse{\pgfmathresult=1}{ 
\tribuPositive{\x*0.75}{0};
}
{
\tribuNegative{\x*0.75}{0};
}
}
\node at (\numberFirstAbs*0.75+0.65,0){$+$};

\foreach \x in {1,...,\numberSecondAbs}  {
\pgfmathsetmacro{\translateX}{(\x+\numberFirstAbs)*0.75+0.5}
\pgfmathparse{ifthenelse(\numberFirst>0,1,0)} \ifthenelse{\pgfmathresult=1}{ 
\tribuPositive{\translateX}{0};
}
{
\tribuNegative{\translateX}{0};
}
}

\node at ( {\numberFirstAbs*0.75+\numberSecondAbs*0.75+1.1},0){$=$};

\pgfmathparse{int(\numberFirstAbs+\numberSecondAbs)}
\foreach \x in {1,...,\pgfmathresult}  {
\pgfmathsetmacro{\translateX}{(\numberFirstAbs+\numberSecondAbs+\x)*0.75+0.95}
\pgfmathparse{ifthenelse(\numberFirst>0,1,0)} \ifthenelse{\pgfmathresult=1}{ 
\tribuPositive{\translateX}{0};
}
{
\tribuNegative{\translateX}{0};
}

}

}
{

\pgfmathparse{int(\numberFirst!=0)}
\ifthenelse{\pgfmathresult = 1}{
\foreach \x in {1,...,\numberFirstAbs}  {

\pgfmathsetmacro{\translateX}{(\x)*0.75}
\pgfmathparse{ifthenelse(\numberFirst>0,1,0)} 
\ifthenelse{\pgfmathresult=1}{
  \tribuPositive{\translateX}{0};
}{
  \tribuNegative{\translateX}{0};
}
}
}{}

\node at (\numberFirstAbs*0.75+0.65,0){$+$};

\pgfmathparse{int(\numberSecond!=0)}
\ifthenelse{\pgfmathresult = 1}{
\foreach \x in {1,...,\numberSecondAbs}  {

\pgfmathsetmacro{\translateX}{(\x+\numberFirstAbs)*0.75+0.5}
\pgfmathparse{ifthenelse(\numberSecond>0,1,0)} 
\ifthenelse{\pgfmathresult=1}{
  \tribuPositive{\translateX}{0};
}{
  \tribuNegative{\translateX}{0};
}
}}{}

\node at ( {\numberFirstAbs*0.75+\numberSecondAbs*0.75+1.1},0){$=$};

\pgfmathparse{ifthenelse(abs(abs(\numberFirst)-abs(\numberSecond))>0,1,0)} 
\ifthenelse{\pgfmathresult=1}{

\node at ( {\numberFirstAbs*0.75+\numberSecondAbs*0.75+1.1},-1){$=$};
\pgfmathsetmacro{\translateX}{(\numberFirstAbs+\numberSecondAbs)*0.75+0.95}

\pgfmathparse{ifthenelse(\numberFirst+\numberSecond>0,1,0)} 
\ifthenelse{\pgfmathresult=1}{
\pgfmathparse{abs(abs(\numberFirst)-abs(\numberSecond))} 
\foreach \x in {1,...,\pgfmathresult}  {
\tribuPositive{\translateX+\x*0.75}{0};
\tribuPositive{\translateX+\x*0.75}{-1};
}
\pgfmathparse{max(abs(\numberFirst),abs(\numberSecond))-abs(abs(\numberFirst)-abs(\numberSecond))} 
\foreach \x in {1,...,\pgfmathresult}  {
\pgfmathsetmacro{\translateX}{(\numberFirstAbs+\numberSecondAbs+abs(abs(\numberFirst)-abs(\numberSecond)))*0.75+0.95}
\explosionForce{\translateX+\x*0.75}{0};
}
}
{
\pgfmathparse{abs(abs(\numberFirst)-abs(\numberSecond))} 
\foreach \x in {1,...,\pgfmathresult}  {
\tribuNegative{\translateX+\x*0.75}{0};
\tribuNegative{\translateX+\x*0.75}{-1};
}
\pgfmathparse{max(abs(\numberFirst),abs(\numberSecond))-abs(abs(\numberFirst)-abs(\numberSecond))} 
\foreach \x in {1,...,\pgfmathresult}  {

\pgfmathsetmacro{\translateX}{(\numberFirstAbs+\numberSecondAbs+abs(abs(\numberFirst)-abs(\numberSecond)))*0.75+0.95}
\explosionForce{\translateX+\x*0.75}{0};
}
}
}
{
\pgfmathsetmacro{\translateX}{(\numberFirstAbs+\numberSecondAbs)*0.75+0.95}
\pgfmathparse{abs(\numberFirst)} 
\foreach \x in {1,...,\pgfmathresult}  {
\explosionForce{\translateX+\x*0.75}{0};
}
}
}
}

\newcommand\equivalentSubtractionAddition[2]{
\pgfmathsetmacro{\numberFirst}{#1}
\pgfmathsetmacro{\numberSecond}{#2}

\pgfmathsetmacro{\numberFirstAbs}{abs(\numberFirst)}
\pgfmathsetmacro{\numberSecondAbs}{abs(\numberSecond)}

\foreach \x in {1,...,\numberFirstAbs}  {
\pgfmathparse{ifthenelse(\numberFirst>0,1,0)} 
\ifthenelse{\pgfmathresult=1}{
\tribuPositive{\x*0.75}{0};
}{
\tribuNegative{\x*0.75}{0};
}
}
\node at (\numberFirstAbs*0.75+0.65,0){$-$};

\foreach \x in {1,...,\numberSecondAbs}  {
\pgfmathparse{ifthenelse(\numberSecond>0,1,0)} 
\ifthenelse{\pgfmathresult=1}{
\tribuPositive{\numberFirstAbs*0.75+\x*0.75+0.5}{0};
}{
\tribuNegative{\numberFirstAbs*0.75+\x*0.75+0.5}{0};
}
}
\node at (\numberFirstAbs*0.75+\numberSecondAbs*0.75+1.1,0){$=$};

\begin{scope}[shift={(\numberFirstAbs*0.75+\numberSecondAbs*0.75+1.1,0)}]
\foreach \x in {1,...,\numberFirstAbs}  {
\pgfmathparse{ifthenelse(\numberFirst>0,1,0)} 
\ifthenelse{\pgfmathresult=1}{
\tribuPositive{\x*0.75}{0};
}{
\tribuNegative{\x*0.75}{0};
}
}
\node at (\numberFirstAbs*0.75+0.65,0){$+$};

\foreach \x in {1,...,\numberSecondAbs}  {
\pgfmathparse{ifthenelse(\numberSecond>0,1,0)} 
\ifthenelse{\pgfmathresult=1}{
\tribuNegative{\numberFirstAbs*0.75+\x*0.75+0.5}{0};
}{
\tribuPositive{\numberFirstAbs*0.75+\x*0.75+0.5}{0};

}
}
\end{scope}
}

\newcommand\subtractionRelativeNumber[2]{
\pgfmathsetmacro{\numberFirst}{#1}
\pgfmathsetmacro{\numberSecond}{#2}

\pgfmathsetmacro{\numberFirstAbs}{abs(\numberFirst)}
\pgfmathsetmacro{\numberSecondAbs}{abs(\numberSecond)}

\foreach \x in {1,...,\numberFirstAbs}  {
\pgfmathparse{ifthenelse(\numberFirst>0,1,0)} 
\ifthenelse{\pgfmathresult=1}{
\tribuPositive{\x*0.75}{0};
}{
\tribuNegative{\x*0.75}{0};
}
}
\node at (\numberFirstAbs*0.75+0.65,0){$-$};

\foreach \x in {1,...,\numberSecondAbs}  {
\pgfmathparse{ifthenelse(\numberSecond>0,1,0)} 
\ifthenelse{\pgfmathresult=1}{
\tribuPositive{\numberFirstAbs*0.75+\x*0.75+0.5}{0};
}{
\tribuNegative{\numberFirstAbs*0.75+\x*0.75+0.5}{0};
}
}
\node at (\numberFirstAbs*0.75+\numberSecondAbs*0.75+1.1,0){$=$};

\begin{scope}[shift={(\numberFirstAbs*0.75+\numberSecondAbs*0.75+1.1,0)}]
\foreach \x in {1,...,\numberFirstAbs}  {
\pgfmathparse{ifthenelse(\numberFirst>0,1,0)} 
\ifthenelse{\pgfmathresult=1}{
\tribuPositive{\x*0.75}{0};
}{
\tribuNegative{\x*0.75}{0};
}
}
\node at (\numberFirstAbs*0.75+0.65,0){$+$};

\foreach \x in {1,...,\numberSecondAbs}  {
\pgfmathparse{ifthenelse(\numberSecond>0,1,0)} 
\ifthenelse{\pgfmathresult=1}{
\tribuNegative{\numberFirstAbs*0.75+\x*0.75+0.5}{0};
}{
\tribuPositive{\numberFirstAbs*0.75+\x*0.75+0.5}{0};

}
}
\end{scope}

\pgfmathsetmacro{\numberAdd}{#1-#2}
\pgfmathsetmacro{\numberAddAbs}{int(abs(\numberAdd)}

\node at (\numberFirstAbs*0.75+\numberSecondAbs*0.75+1.1,-1){$=$};

\foreach \x in {1,...,\numberAddAbs}  {
\pgfmathsetmacro{\position}{\x*0.75+\numberFirstAbs*0.75+\numberSecondAbs*0.75}

\pgfmathparse{ifthenelse(\numberAdd>0,1,0)} 
\ifthenelse{\pgfmathresult=1}{
\tribuPositive{\position+1.1}{-1};
}{
\tribuNegative{\position+1.1}{-1};
}
}

}

\title{Two teams that cancel:\\ a didactic analysis of a teaching progression\\ for signed numbers from Grade~5 to Grade~7}
\author{Vincent Tariel-Adam\\[2pt]
  \small Independent researcher; CommeUnJeu (\url{https://commeunjeu.org}), France\\
  \small \texttt{vincent.tariel@alumni.polytechnique.org}}
\date{October 2026}

\begin{document}
\maketitle

\begin{abstract}
\noindent
Signed numbers make learners revise what addition, subtraction and the minus sign meant for whole numbers. This paper presents the design rationale and a critical analysis of a three-level source progression for signed numbers (Grades 5--7, nominal ages 10--13), from which French-curriculum and International Baccalaureate chapters published at \url{https://commeunjeu.org} are assembled. One representation is presented first and taken up throughout: an integer is drawn as a collection of positive and negative units, a positive and a negative unit forming a zero pair that leaves nothing. An integer therefore has many drawings, which makes visible the construction of the integers as classes of pairs of natural numbers. The design hypothesis is that the gesture pupils perform with the units, inserting a zero in the form $a+(-a)$, is shared by the subtraction of a negative number, by the worked examples that accompany the 2026 French programme, and by the derivation of the sign of a product of two negatives through distributivity; the progression does not model that product, it derives it from the zero pair, a derivation that holds only if distributivity is kept. The number line supports position, order and distance; contexts provide applications. The choices are related to the epistemological obstacles described by Glaeser, to studies of counter and number-line models, including a classroom comparison in which a path model gave better results than chips, immediately and five weeks later, and to the French and IB curricula. A block-level audit of the dated source documents found the representation absent from the first level, no exercise asking pupils to draw, and the decisive subtraction never drawn; the revisions made after the audit are recorded. No empirical evaluation is reported; hypotheses for one are stated.
\end{abstract}


\section{Introduction}\label{sec:intro}

Signed numbers are a turning point of school arithmetic. With them, a subtraction is always possible, adding can make a number smaller and subtracting can make it larger, the order of the magnitudes reverses below zero ($-5<-2$ although $5>2$), and the minus sign takes three roles: the operation of subtraction, the sign of a number, and the opposite of a number \citep{Vlassis2004,Vlassis2008,Lamb2012}. According to the abstract of Glaeser's study, mathematicians needed more than fifteen centuries before the rule of signs seemed trivial \citep{Glaeser1981}. Research on learners is comparatively thin: a review found 47 research-journal publications on integers between 1980 and 2012 \citep{Bishop2014}.

A curriculum author meets a design problem specific to this domain. The first representation chosen for negative numbers must carry their existence, their order, addition and subtraction, and it must not stand in the way of multiplication. The didactic literature diagnoses a recurring difficulty at that last point: concrete additive models (gains and debts, rising and falling temperatures, moves along a line) serve addition, but extending them to the product of two negatives requires additional conventions and may create difficulties, up to blocking the product \citep{Glaeser1981,IREMAquitaine2008,Eduscol2016Relatifs,Matheron2016,Gliere2025}.

This paper documents and analyses one answer to that problem, a three-level progression written as bilingual (English and French) lessons and exercise sets. One representation is presented first and taken up throughout: an integer is a collection of units of two kinds, positive and negative, and one positive unit with one negative unit forms a \emph{zero pair} that leaves nothing. The paper's thesis is not that this representation shows everything; it is that it makes one gesture available to the pupil, \emph{inserting a zero in the form $a+(-a)$}, and that the same gesture underlies the subtraction of a negative number and the derivation of the sign of a product of two negatives by distributivity (Figure~\ref{fig:zero}). The design rests on five principles:
\begin{enumerate}[label=(\roman*)]
  \item an integer has several drawings, and every balanced collection represents zero; the opposite is the same units of the other kind (Section~\ref{sec:p1});
  \item adding and subtracting are gestures on the units, and subtracting a number that is not there begins by inserting a zero (Section~\ref{sec:p2});
  \item the sign of a product of two negatives is derived from the zero pair by distributivity, not read on a model of the product (Section~\ref{sec:p3});
  \item each representation is kept for what it carries: situations for the encounter, the number line for position, order and distance, the units for equivalence and compensation; and the drawing is abandoned when the rules can stand alone (Section~\ref{sec:p4});
  \item the drawing comes first, then the chain of equalities, written line for line under it (Section~\ref{sec:p5}).
\end{enumerate}

\paragraph{Genre and scope.} The paper is a design rationale and a curriculum analysis, the companion of a study of a fraction progression by the same author \citep{TarielAdam2026Fractions}, whose method it shares. Like a hypothetical learning trajectory \citep{simon1995,clements2004}, the progression embodies learning goals, an order of activities and hypotheses about learning; the paper makes them explicit and confronts them with the literature and with the documents that implement them. Its method borrows the \emph{a priori analysis} of didactic engineering \citep{artigue1990,artigue1992} at the level of global choices (the representation, the order of meanings) and, in Section~\ref{sec:task}, at the level of one task and its didactic variables. No classroom experiment is reported; statements about learners are taken from the cited literature or presented as design hypotheses. None of the five principles is new taken alone: the counter or ``charge'' model is classical \citep{Battista1983}, and the derivation of the sign rule by distributivity is the route of the French official documents \citep{Eduscol2016Relatifs,MEN2026Exemples}. What the paper contributes is the articulation of one representation across three levels, the identification of the gesture that links its additive and multiplicative uses, the explicit statement of where the drawing stops, and a confrontation of the stated design with the materials, which led to revisions made before this paper was written (Section~\ref{sec:revisions}).

\paragraph{Grade convention.} The progression is organised in three source levels, labelled Grade~5 to Grade~7 and nominally corresponding to ages 10 to 13. It is a source progression, not a national curriculum: curriculum chapters are assembled from it by links (Section~\ref{sec:method}). In the French-curriculum chapters, Grade~5 and Grade~6 are both used in \emph{cinqui\`eme} (age 12--13) and Grade~7 in \emph{quatri\`eme}; in the International Baccalaureate chapters, Grade~5 is used in PYP~5, Grade~6 in MYP~1 and Grade~7 in MYP~2 (Table~\ref{tab:projections}). Unless stated otherwise, ``Grade~$N$'' refers to the levels of the progression.

\section{Corpus and method}\label{sec:method}

\paragraph{Corpus.} The corpus is the six source documents of the signed-number progression: a lesson and an exercise set for each of Grades 5, 6 and 7. They are written in a mark-up in which every definition, method, example, exercise and answer is a delimited \emph{block} with a type and a stable identifier; a block of a curriculum chapter, or of a higher level, may be a \emph{link} to a block written elsewhere, which it displays without copying. Two dated states are analysed: the state of 28~October~2026 before the revisions, on which the audit was made, and the state after the revisions of the same day (Section~\ref{sec:revisions}). Table~\ref{tab:corpus} gives their sizes.

\paragraph{Statements and counts.} As in the companion paper, three kinds of statement are made about the corpus: \emph{presence} (a block with a given content exists at a given level), \emph{order} (one block precedes another), and \emph{counts}, exhaustive over the mark-up. Counts are made by a script on the English half of each block and checked by reading; the French half is a parallel translation. A block \emph{draws units} when it uses the drawing macros of the unit model. A block \emph{cites a situation} when its English text names one of a closed list of contexts (temperatures, heights and depths, floors, scores, money); the list and the blocks it selects are kept with the script. A \emph{chain} is a multi-line calculation; it carries a \emph{side annotation} when at least one line has a note in the margin, a \emph{verbal reason} when that note contains words rather than only a side calculation, and it is \emph{fully annotated} when every line but the last carries a note. An annotation is not a justification: whether a note gives a valid reason is decided by reading, not by the script. Content counts inspect the literal body of each block without resolving links: a linked copy whose body is filled in is counted, an empty link contributes no content; a fall in a count can therefore come from a block that became a link, not from a change in what a chapter displays. The two states are identified by two commits of the author's repository (\texttt{3a06d7d472} and \texttt{0ce41eca1c}).

\paragraph{Projections.} The counts concern the source documents, not repeated appearances of a block in the assembled chapters. The curriculum chapters (French \emph{cinqui\`eme} and \emph{quatri\`eme}; IB PYP~5, MYP~1 and MYP~2) are described from their link tables, which list, unit by unit, the blocks they show; they establish what a chapter contains and in which order, not what pupils actually experienced.

\paragraph{Limits of a self-audit.} The author designed the progression and made the audit. Two external readings of the revision proposals and one of the plan of this paper were made by a language-model reviewer, which found errors in the author's first coding (a regular expression had missed one of the two written forms of a margin note, and a word list had matched ``owe'' inside ``lower''); the corrected counts are those reported. The source documents and the counting script can be consulted on request; the assembled chapters are public (Appendix~\ref{app:resources}). No public archive of the dated sources accompanies this manuscript, so the counts are not yet independently reproducible from the public chapters alone.

\begin{table}[htbp]
\centering\small
\caption{The six source documents before and after the revisions of 28~October~2026 (English halves). ``Units'' counts blocks that draw positive or negative units. Exercise counts exclude multiple-choice questions (4, 12 and 9 at Grades 5, 6 and 7); links are blocks that display a block written elsewhere.}
\label{tab:corpus}
\begin{tabular}{lrrr}
\toprule
 & Grade 5 & Grade 6 & Grade 7 \\
\midrule
Lesson blocks, before $\to$ after & $18\to21$ & $57\to57$ & $81\to81$ \\
\quad of which links, after & 0 & 3 & 21 \\
Lesson blocks drawing units, before $\to$ after & $0\to5$ & $7\to6$ & $10\to11$ \\
Lesson blocks citing a situation, before & 15 & 3 & 1 \\
Exercises, before $\to$ after & $8\to27$ & $171\to179$ & $219\to219$ \\
\quad of which links, after & 0 & 15 & 109 \\
Written exercises (photo or slate), before $\to$ after & $4\to8$ & $33\to41$ & $38\to38$ \\
\quad of which ask the pupil to draw units, after & 4 & 8 & 0 \\
Multi-line chains in exercise answers, before $\to$ after & $0\to0$ & $68\to73$ & $94\to94$ \\
\quad with a side annotation, before $\to$ after & -- & $28\to33$ & $48\to48$ \\
\quad with a verbal reason, before $\to$ after & -- & $21\to26$ & $41\to41$ \\
\bottomrule
\end{tabular}
\end{table}

\section{Theoretical background}\label{sec:background}

\subsection{Obstacles}
Glaeser's study of the history of the rule of signs, from Diophantus to the nineteenth century, identifies a series of obstacles \citep{Glaeser1981}. His text is not openly available to the author; \citet{SalvadorNacarato2003} reproduce six obstacles, which \citet{Bishop2014} corroborate and \citet{IREMAquitaine2008} regroup into four: the inability to manipulate isolated negative quantities; the difficulty of giving them a meaning; the difficulty of unifying the number line; the ambiguity of the two zeros, absolute and origin; stagnation at the stage of concrete operations; and the desire for a unifying model, that is, the use of an additive model in the multiplicative field, where it does not apply. A French classroom study built on Glaeser reports a pupil who maintained $(-3)\times(-5)=-15$ because descending five steps three times lands at $-15$, and concludes that a strong ``going up and down'' image becomes an obstacle to multiplication \citep{IREMAquitaine2008}. In history and in children's thinking, \citet{Bishop2014} find the same obstacles (quantities less than nothing, removing more than one has, the belief that addition cannot make smaller) and also affordances: order, magnitude, computational reasoning, and logical necessity. \citet{Gallardo2002} extracts from historical texts four levels of acceptance of negative numbers (as subtrahend, relative number, isolated number, formal number) and observed the first three among pupils aged 12 to 13; \citet{CoquinViennot1985} describes an ordered hierarchy of conceptions of the integers among French pupils of the four years of \emph{coll\`ege}. Progress through these conceptions is slow, and a concrete representation has to be designed as a passage, not as a destination.

\subsection{The minus sign}
\citet{Vlassis2008} attributes many difficulties with negatives in equations less to the concept than to a restricted use of the minus sign, and \citet{Vlassis2004} describes the flexibility in ``negativity'' that pupils must develop to use the sign in its several roles. Explicit teaching of the unary and binary meanings of ``$-$'' changed the mental models of first-graders \citep{Bofferding2014}, and teachers can support each interpretation \citep{Lamb2012}. A representation with two labelled kinds of units writes the sign of a number on the object itself; the passage to written equalities is where its roles have to be negotiated.

\subsection{Models and their breaking points}
\citet{Janvier1985} sorts the models of textbooks into two families: number-line models, where an integer is a position or a motion, and models that combine objects of opposite nature and represent an integer by pairs, that is, by an equivalence class. In his teaching experiment with pupils aged 13--14, an image model of cancellation (a balloon with sand bags and helium balloons) did significantly better than a symbolic model on addition, with no difference on subtraction, which in such models, he notes, often requires ``acrobatics''. \citet{Battista1983} presented the positive-negative charge model as concrete and complete for the four operations. \citet{Liebeck1990} taught a model of scores and forfeits and compared it with a number-line model; \citet{LinchevskiWilliams1999} used situations from everyday life to fill the gap between operational and structural conceptions. \citet{Vig2014} argue that every model has a breaking point and analyse the chip model of integer subtraction: when too few chips are present, zero pairs must be added, and the use of the model can turn into the execution of a graphical procedure that the pupil does not understand.

The most direct comparison found is \citet{NurnbergerHaag2015PMENA}. Eight intact classes of first-time learners in the United States (fifth and sixth grade; 154 pupils analysed) were randomly assigned to a chip model or to a ``walk-it-off'' number-line model, for eight lessons on the four operations. Both groups learned; the number-line group scored about 11 points higher than the chip group on a 100-point test at the post-test and about 13 points higher five weeks later. The difference concerned a pooled set of 14 problems in which the chips require movements inconsistent with the operation (subtractions and multiplications such as $-2-(-5)$ and $-4\times(-3)$, where zero pairs must be inserted before chips can be removed); on the 22 other problems no significant difference was found. The chip group showed greater knowledge of sums of opposites at the post-test, a difference that was no longer significant five weeks later; this knowledge was assessed with calculations, generated problems and explanations of why, for instance, $-8+(-7+7)$ and $-8+(-5+5)$ are equal. For the 118 of 154 pupils without strong prior knowledge of opposite sums, the chip group scored about 13 points lower at the post-test and about 16 points lower five weeks later. The author also notes that chips represent zero by several objects, against the expectation of ``no thing'' \citep{NurnbergerHaag2015PMENA}. These findings concern a particular chip model and a particular path model; they are discussed in Sections~\ref{sec:task} and~\ref{sec:discussion}. Other studies support number-line routes: a multimedia number-line metaphor reduced the confusion of the negative sign with the subtraction operator \citep{MorenoMayer1999}, and a design experiment built addition and subtraction on a financial context and an empty vertical line \citep{StephanAkyuz2012}. The order relation, which a collection of units does not show, is studied by \citet{Schindler2017}.

\subsection{The permanence of formal laws}
\citet{Hankel1867} introduced the negative integers as formal differences and stated that the products $(-a)c=-ac$ and $(-a)(-c)=ac$ cannot be proved in formal mathematics: they are conventions adopted to preserve the laws of calculation, chief among them distributivity, from which he derives the sign rule. According to its abstract, \citet{HefendehlHebeker1991} views negative numbers as an extension of the number system that overcomes the historical obstacles; \citet{Bishop2014} discuss the permanence principle with reference to Hefendehl-Hebeker and to Freudenthal's didactical phenomenology \citep{Freudenthal1983}, and count this ``logical necessity'' among the productive resources of pupils.

\subsection{States and transformations}
\citet{Vergnaud1982} distinguishes, in additive problems, measures, transformations and static relationships; directed numbers represent the latter two. In his terms a temperature below zero is an abscissa on an axis rather than a measure, and a fall of four degrees is a transformation; reading the first as the second is the confusion at stake. A representation by units has to say what it draws.

\subsection{Institutional context}
In France, a new programme for \emph{cycle~4} was published in March 2026 \citep{MEN2026Arrete}; according to its published calendar it applies to \emph{cinqui\`eme} from 2026--2027 and to \emph{quatri\`eme} from 2027--2028. The examples that accompany it \citep{MEN2026Exemples} define the opposite of a decimal $a$ as the unique number $b$ such that $a+b=0$; establish the rule of addition from the cancellation of opposites, for instance $6+(-9)=6-6+(-3)=-3$; establish that subtracting is adding the opposite by inserting a zero, as in $3{,}1-(-2)=3{,}1+2+(-2)-(-2)=5{,}1$; obtain the sign of a product of two negatives by distributivity over a sum equal to zero; and keep the number line for order and distance to zero. The resources of 2016 for the previous programme \citep{Eduscol2016Relatifs} proposed to introduce negative numbers without contexts, from the equivalence of calculations such as $3-5=2-4=1-3=0-2$ whose common result is written $-2$, warned that the models of gains and losses or of rising and falling ``can be a real obstacle'' to multiplying two negatives, advised against relying on a single model such as temperature, and justified the sign rule by distributivity, provided that distributivity is extended. Two French groups made similar diagnoses and proposed routes through equalities with a missing term \citep{IREMAquitaine2008} or through calculation programmes \citep{Matheron2016}. In the International Baccalaureate, the Middle Years Programme prescribes no textbook; the publisher's sequence examined here places the addition and subtraction of negative numbers on number lines, and pattern-based multiplication, in MYP~1 \citep{Haese2022Scope}.

\section{Principle 1: an integer has several drawings, every balanced collection is zero}\label{sec:p1}

\begin{principle}
An integer is drawn as a collection of positive and negative units. A positive and a negative unit form a zero pair, which leaves nothing; the integer is what is left, with its sign, once the zero pairs are set aside. An integer therefore has many drawings, every collection with as many positive as negative units represents zero, and the opposite of an integer is the same units of the other kind.
\end{principle}

\subsection{Three entry points}
Three entries into signed numbers are documented in the materials and the literature reviewed in Section~\ref{sec:background}; Table~\ref{tab:entries} summarises what each makes immediately available. (a)~\emph{Situations and the number line}: a negative number is a level below a zero (a temperature, a depth, a floor), placed on a line drawn upright or across; this is the entry of many textbooks and of the IB scope and sequence cited above. (b)~\emph{Differences}: a negative number is the common result of calculations such as $3-5=2-4=1-3=0-2$, or the missing term of $5+\ldots=2$; this is the formal route of the 2016 French resources and of Hankel. (c)~\emph{Units that cancel}: a negative number is a collection of negative units, and a positive and a negative unit cancel \citep{Janvier1985,Battista1983,Liebeck1990}.

\begin{table}[htbp]
\centering\small
\caption{What each entry makes immediately available, as described in the cited sources and in the design; ``--'' means that the entry needs an additional convention, not that it is impossible.}
\label{tab:entries}
\begin{tabularx}{\textwidth}{>{\raggedright\arraybackslash}p{3.0cm}XXX}
\toprule
 & (a) situations and line & (b) differences & (c) units that cancel \\
\midrule
Encounter, need & yes: a level below zero & -- & -- \\
Order & yes: lower or further left is smaller & -- & -- (a collection shows no order) \\
Opposite & symmetry about zero & $a+b=0$ & the same units of the other kind; together they leave nothing \\
Several writings of one number & -- & yes: $3-5=0-2$ & yes: several drawings \\
Addition & moves (positions and motions mixed) & formal & put together, cancel pairs \\
Subtraction & reversing a move & adding the opposite & take away; insert zero pairs if needed \\
Product of two negatives & requires additional conventions & permanence, distributivity & not modelled in this progression; justified by distributivity from the zero pair \\
\bottomrule
\end{tabularx}
\end{table}

The three entries are not mathematically opposed. A drawing with $p$ positive and $m$ negative units codes the pair of natural numbers $(p,m)$; two drawings represent the same integer exactly when $p+m'=p'+m$, which is the equivalence that defines the integers as classes of pairs, the construction implicit in Hankel's formal differences and in the 2016 sequence $3-5=2-4=\ldots=0-2$. Putting two drawings together adds the pairs, $(p,m)+(p',m')=(p+p',m+m')$; exchanging the two kinds of units gives the opposite, $(m,p)$; adding a zero pair, $(p+1,m+1)$, does not change the class. Entry~(c) is the drawable form of entry~(b). The identification has limits that the paper does not hide. The order and the product can be defined on pairs ($(p,m)\le(p',m')$ when $p+m'\le p'+m$; $(p,m)(q,n)=(pq+mn,pn+mq)$), but they require rules that the gesture of cancelling does not give; and units worth one tenth would draw decimals, so that stopping the drawing at the integers is a didactic choice, not a mathematical necessity (Section~\ref{sec:p4}).

\begin{figure}[htbp]
\centering
\begin{minipage}{0.45\textwidth}\centering
\begin{tikzpicture}
            \node at (-1,0) {$\textcolor{colordef}{+2}=$};
            \tribuPositive{0}{0};
            \tribuPositive{0.75}{0};
            \end{tikzpicture}

\medskip
\begin{tikzpicture}
            \node at (-1,0) {$\textcolor{colorprop}{-3}=$};
            \tribuNegative{0}{0};
            \tribuNegative{0.75}{0};
            \tribuNegative{1.5}{0};
            \end{tikzpicture}
\end{minipage}\hfill
\begin{minipage}{0.52\textwidth}\centering
\begin{tikzpicture}
        \tribuPositive{-1.5}{0}; \tribuPositive{-0.75}{0}; \tribuPositive{0}{0};
        \node at (-0.75,-2) {$\textcolor{colordef}{(+3)}$};
        \node at (0.75,0) {$+$}; \node at (0.75,-2) {$+$};
        \tribuNegative{1.5}{0};
        \node at (1.5,-2) {$\textcolor{colorprop}{(-1)}$};
        \node at (2.25,0) {$=$}; \node at (2.25,-1) {$=$}; \node at (2.25,-2) {$=$};
        \tribuPositive{3}{0}; \tribuPositive{3.75}{0}; \explosionForce{4.5}{0};
        \tribuPositive{3}{-1}; \tribuPositive{3.75}{-1};
        \node at (3.37,-2) {$\textcolor{colordef}{+2}$};
        \end{tikzpicture}

\medskip
\begin{tikzpicture}
        \tribuPositive{-0.75}{0}; \tribuPositive{0}{0};
        \node at (-0.5,-1) {$\textcolor{colordef}{(+2)}$};
        \node at (0.75,0) {$+$}; \node at (0.75,-1) {$+$};
        \tribuNegative{1.5}{0}; \tribuNegative{2.25}{0};
        \node at (1.75,-1) {$\textcolor{colorprop}{(-2)}$};
        \node at (3,0) {$=$}; \node at (3,-1) {$=$};
        \explosionForce{3.75}{0}; \explosionForce{4.5}{0};
        \node at (4.25,-1) {$0$};
        \end{tikzpicture}
\end{minipage}
\caption{The representation as the Grade~5 lesson draws it (opening problem, block 77706). Left: the sign in front of a number says which kind of unit a group contains. Right: ``when 3 positives meet 1 negative'' and ``when 2 positives meet 2 negatives''; a cancelled zero pair is drawn as a burst, each line is an equality between drawings, and the numbers are written underneath, under the drawings they translate.}
\label{fig:particles}
\end{figure}
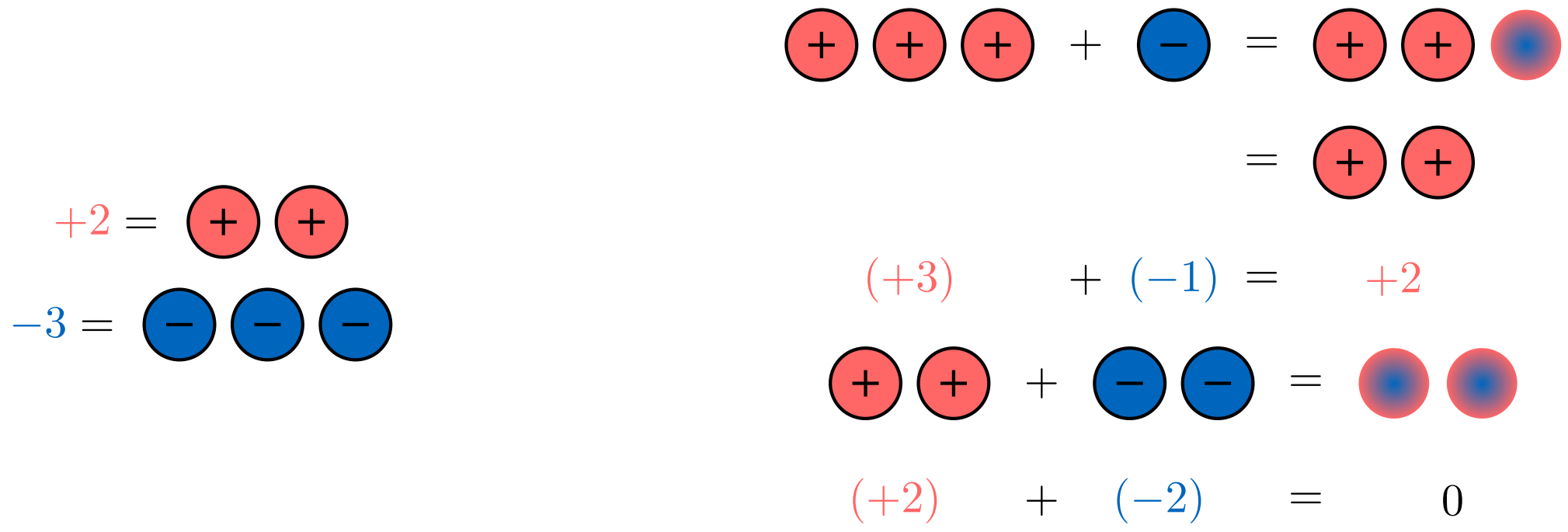

\subsection{How the progression enters}
The revised Grade~5 lesson opens with this representation, before any situation (Figure~\ref{fig:particles}) (Section~\ref{sec:revisions} gives the earlier state). An opening problem imagines ``a world with two types of particles'', positives and negatives, each standing for one unit: particles of the same type join forces, a positive and a negative particle cancel each other, and this is a zero pair. The sign written in front of a number says which type of particle a group contains. Each case is drawn as an equality between drawings, with the result in numbers underneath ($(+3)+(-1)=+2$; $(+2)+(-2)=0$). A definition follows: the integers are the whole numbers and their opposites; positive numbers carry a plus sign, often omitted; two integers are opposites if, drawn without zero pairs, they have the same number of units with opposite signs; placing $(+2)$ and $(-2)$ together leaves no unit. Zero is described according to the convention of each language: neither positive nor negative in English, both positive and negative in French, where it is also neither strictly positive nor strictly negative.

Only then come the situations. A thermometer, a diver and a car park ask how to write three positions with one number each; a method, ``reading a signed number in a situation'', writes a quantity below zero with a negative number, and its table adds to the rows \emph{temperature}, \emph{height} and \emph{building} a row \emph{units}, where $3$ above zero is three positive units and $3$ below zero three negative units. The number line follows, upright like a thermometer or across like a ruler, with the remark that, with a step of~1, the number of steps from zero is the number of units. The definition of opposite numbers on the line, which preceded the units in the earlier state, has become a remark on the line, which begins with the units: two opposite numbers, drawn without zero pairs, have the same number of units with opposite signs, and $(+4)+(-4)=0$.

\subsection{What the choice gains and what it postpones}
The choice gives the pupil from the first lesson an object to which the opposite, zero, and the equality of different writings are attached, and a gesture (Principle~2) that will be reused. It postpones two things. The \emph{order} is not visible on a collection: $-5<-2$ is read on the line, by a separate rule (Principle~4). The \emph{situations} are no longer the definition: pupils are asked to coordinate a position relative to a chosen zero ($-4\,^\circ$C, four degrees below zero) with a collection of signed units (four negative units); whether pupils of ten do so without reading a position as a change is an open question \citep{Vergnaud1982}. The choice also has a cost that the literature names: a collection with several objects represents nothing, against the expectation that zero is ``no thing'' \citep{NurnbergerHaag2015PMENA}. The design accepts that cost deliberately, because the several drawings of zero are what Principle~2 uses; it is a hypothesis, not a demonstrated gain.

\section{Principle 2: inserting a zero}\label{sec:p2}

\begin{principle}
Adding is putting two collections together and setting the zero pairs aside. Subtracting is taking units away; when the units to take away are not there, zero pairs are first inserted, which does not change the number. The rules of addition and subtraction are read on the drawings and then stated with the distance to zero, so that they survive when the drawing stops.
\end{principle}

\subsection{Adding}
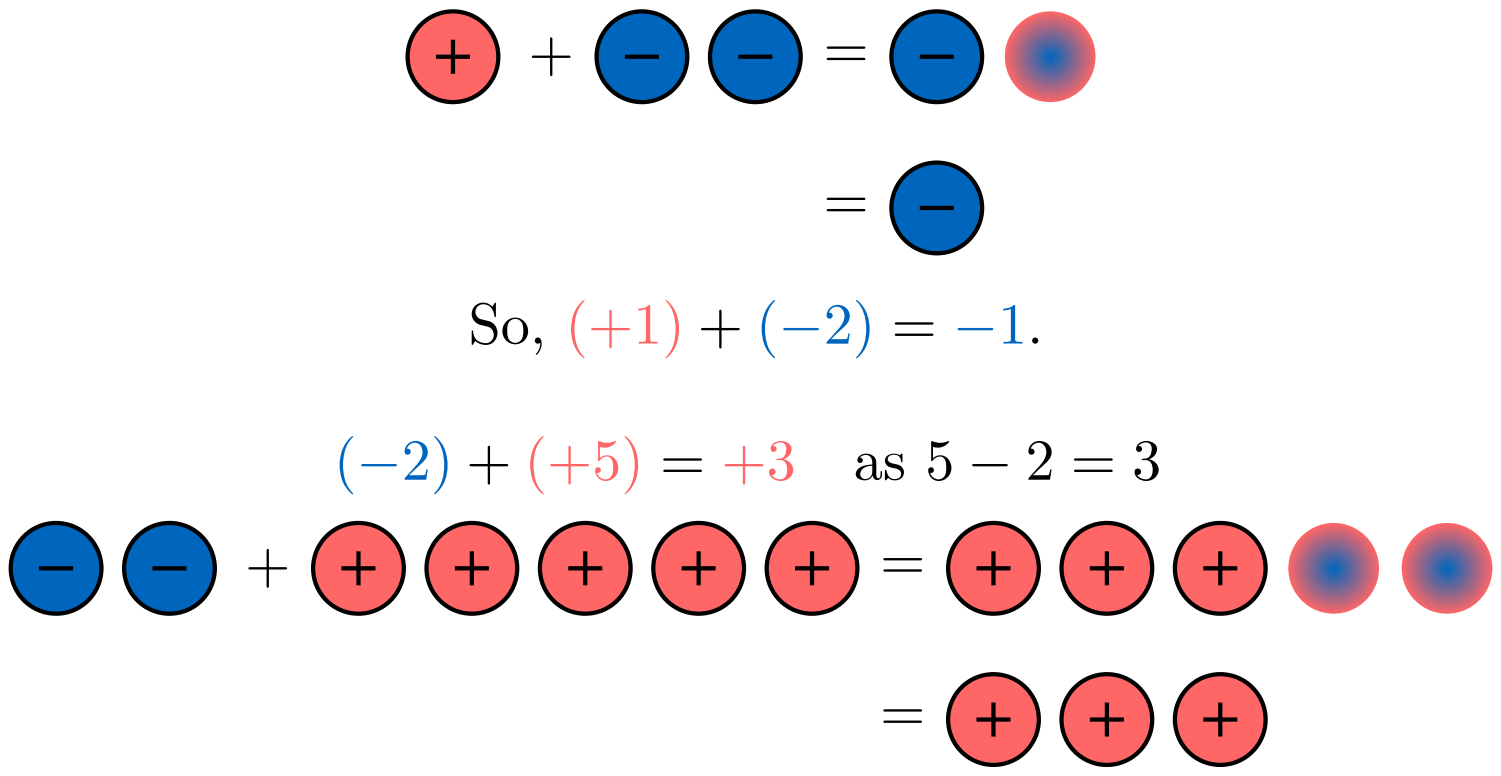
\begin{figure}[htbp]
\centering
\begin{tikzpicture}\additionRelativeNumber{1}{-2}\end{tikzpicture}

\smallskip
So, $\textcolor{colordef}{(+1)}+\textcolor{colorprop}{(-2)}=\textcolor{colorprop}{-1}$.

\bigskip
$\textcolor{colorprop}{(-2)}+\textcolor{colordef}{(+5)}=\textcolor{colordef}{+3} \quad \text{as }5-2=3$\\[4pt]
\begin{tikzpicture}
\additionRelativeNumber{-2}{5}
\end{tikzpicture}
\caption{Additions as the Grade~6 lesson draws them (worked example 77708--77709 and method 77711): the units of the two numbers, then the same units with the zero pairs cancelled, then what is left; the result in numbers below or above.}
\label{fig:addition}
\end{figure}

The Grade~6 lesson draws each addition as an equality between drawings: the units of the first number, the units of the second, an equal sign, the same units with the zero pairs replaced by small bursts, and, on the next line, what is left (Figure~\ref{fig:addition}). The method then states the rules with the distance to zero: two positives or two negatives, add the distances and keep the sign; a positive and a negative, subtract the smaller distance from the larger and take the sign of the number with the larger distance; equal distances give zero. A proposition, ``regrouping by sign'', applies the two teams at the level of a calculation: in $(+16)+(-11)+(+4)+(-9)$, the positives are gathered, then the negatives, and $(+20)+(-20)=0$. A remark warns that adding can make smaller: it is the sign of the number added, not the operation, that decides.

\begin{figure}[htbp]
\centering
\begin{minipage}{0.48\textwidth}\centering
\small $\textcolor{colordef}{(+3)}-\textcolor{colordef}{(+2)}$: we remove $2$ positives from $3$ positives.\\[4pt]
\begin{tikzpicture}
\pgfmathsetmacro{\numberFirst}{3}
\pgfmathsetmacro{\numberSecond}{2}
\foreach \x in {1,...,\numberFirst}  {
\tribuPositive{\x*0.75}{0};
}
\node at (\numberFirst*0.75+0.65,0){$-$};
\foreach \x in {1,...,\numberSecond}  {
\tribuPositive{\numberFirst*0.75+\x*0.75+0.5}{0};
}
\node at (\numberFirst*0.75+\numberSecond*0.75+1.1,-1){$=$};
\node at (\numberFirst*0.75+\numberSecond*0.75+1.1,0){$=$};
\pgfmathsetmacro{\numberdiff}{int(\numberFirst-\numberSecond)}
\foreach \x in {1,...,\numberdiff}  {
\tribuPositive{\numberFirst*0.75+\numberSecond*0.75+\x*0.75+1}{0};
\tribuPositive{\numberFirst*0.75+\numberSecond*0.75+\x*0.75+1}{-1};
}
\foreach \x in {1,...,\numberSecond}  {
\pgfmathsetmacro{\translate}{\numberFirst*0.75+\numberSecond*0.75+\numberdiff*0.75+\x*0.75+1}
\tribuPositive{\translate}{0};
\draw[ultra thick] (\translate-0.4,-0.4)--(\translate+0.4,0.4);
\draw[ultra thick] (\translate+0.4,-0.4)--(\translate-0.4,0.4);
}
\end{tikzpicture}
\end{minipage}\hfill
\begin{minipage}{0.48\textwidth}\centering
\small $\textcolor{colordef}{(+3)}+\textcolor{colorprop}{(-2)}$: we again remove $2$ positives.\\[4pt]
\begin{tikzpicture}
\additionRelativeNumber{3}{-2};
\end{tikzpicture}
\end{minipage}

\bigskip
\small $\textcolor{colordef}{(+3)}-\textcolor{colordef}{(+2)}=\textcolor{colordef}{(+3)}+\textcolor{colorprop}{(-2)}$\\[4pt]
\begin{tikzpicture}
\equivalentSubtractionAddition{3}{2};
\end{tikzpicture}
\caption{The opening problem of subtraction (block 77714, Grade~6), first case: units taken away are crossed out, a cancelled pair is a burst, and the two operations leave the same units.}
\label{fig:takeaway}
\end{figure}
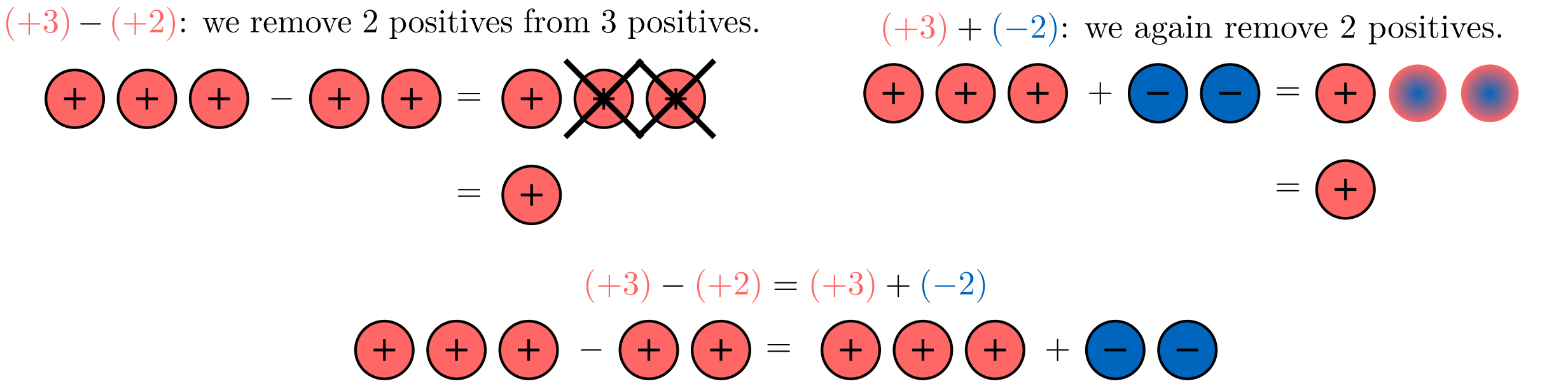

\subsection{Subtracting, and the decisive case}\label{sec:subtraction}
\begin{figure}[htbp]
\centering
\small Before the revision (model answer 77719):\\[2pt]
$$\begin{aligned}[t]
\textcolor{colordef}{(+4)}-\textcolor{colorprop}{(-2)}&=\textcolor{colordef}{(+4)}+\textcolor{colordef}{(+2)}&&\text{(add the opposite)}\\
&=\textcolor{colordef}{+6}&&\text{(same sign: add the distances to zero)}
\end{aligned}$$
\begin{tikzpicture}\subtractionRelativeNumber{4}{-2};\end{tikzpicture}

\bigskip
\small After the revision (block 77714, third case):\\[4pt]
\begin{tikzpicture}[scale=0.9]
\tribuPositive{0.75}{0.00}
\tribuPositive{1.50}{0.00}
\tribuPositive{2.25}{0.00}
\tribuPositive{3.00}{0.00}
\node at (3.65,0.00){$-$};
\tribuNegative{4.60}{0.00}
\tribuNegative{5.35}{0.00}
\node at (5.95,0.00){$=$};
\tribuPositive{6.80}{0.00}
\tribuPositive{7.55}{0.00}
\tribuPositive{8.30}{0.00}
\tribuPositive{9.05}{0.00}
\tribuPositive{9.80}{0.36}
\tribuNegative{9.80}{-0.36}
\tribuPositive{10.55}{0.36}
\tribuNegative{10.55}{-0.36}
\draw[dashed, rounded corners] (9.38,-0.80) rectangle (10.97,0.80);
\node at (11.20,0.00){$-$};
\tribuNegative{12.15}{0.00}
\tribuNegative{12.90}{0.00}
\node at (5.95,-2.00){$=$};
\tribuPositive{6.80}{-2.00}
\tribuPositive{7.55}{-2.00}
\tribuPositive{8.30}{-2.00}
\tribuPositive{9.05}{-2.00}
\tribuPositive{9.80}{-1.64}
\tribuNegative{9.80}{-2.36}
\tribuPositive{10.55}{-1.64}
\tribuNegative{10.55}{-2.36}
\draw[dashed, rounded corners] (9.38,-2.80) rectangle (10.97,-1.20);
\draw[ultra thick] (9.40,-2.76)--(10.20,-1.96);
\draw[ultra thick] (10.20,-2.76)--(9.40,-1.96);
\draw[ultra thick] (10.15,-2.76)--(10.95,-1.96);
\draw[ultra thick] (10.95,-2.76)--(10.15,-1.96);
\node at (5.95,-3.60){$=$};
\tribuPositive{6.80}{-3.60}
\tribuPositive{7.55}{-3.60}
\tribuPositive{8.30}{-3.60}
\tribuPositive{9.05}{-3.60}
\tribuPositive{9.80}{-3.60}
\tribuPositive{10.55}{-3.60}
\end{tikzpicture}

\caption{The decisive subtraction, $(+4)-(-2)$. Above, the earlier model answer: the drawing already applies the rule (``add the opposite'') that it should justify. Below, the case added to the opening problem by the revision: there is no negative unit to take away, so two zero pairs are inserted (framed), the first line still equals the starting number minus $(-2)$, the two negative units are then crossed out, and six positive units are left.}
\label{fig:subtraction}
\end{figure}
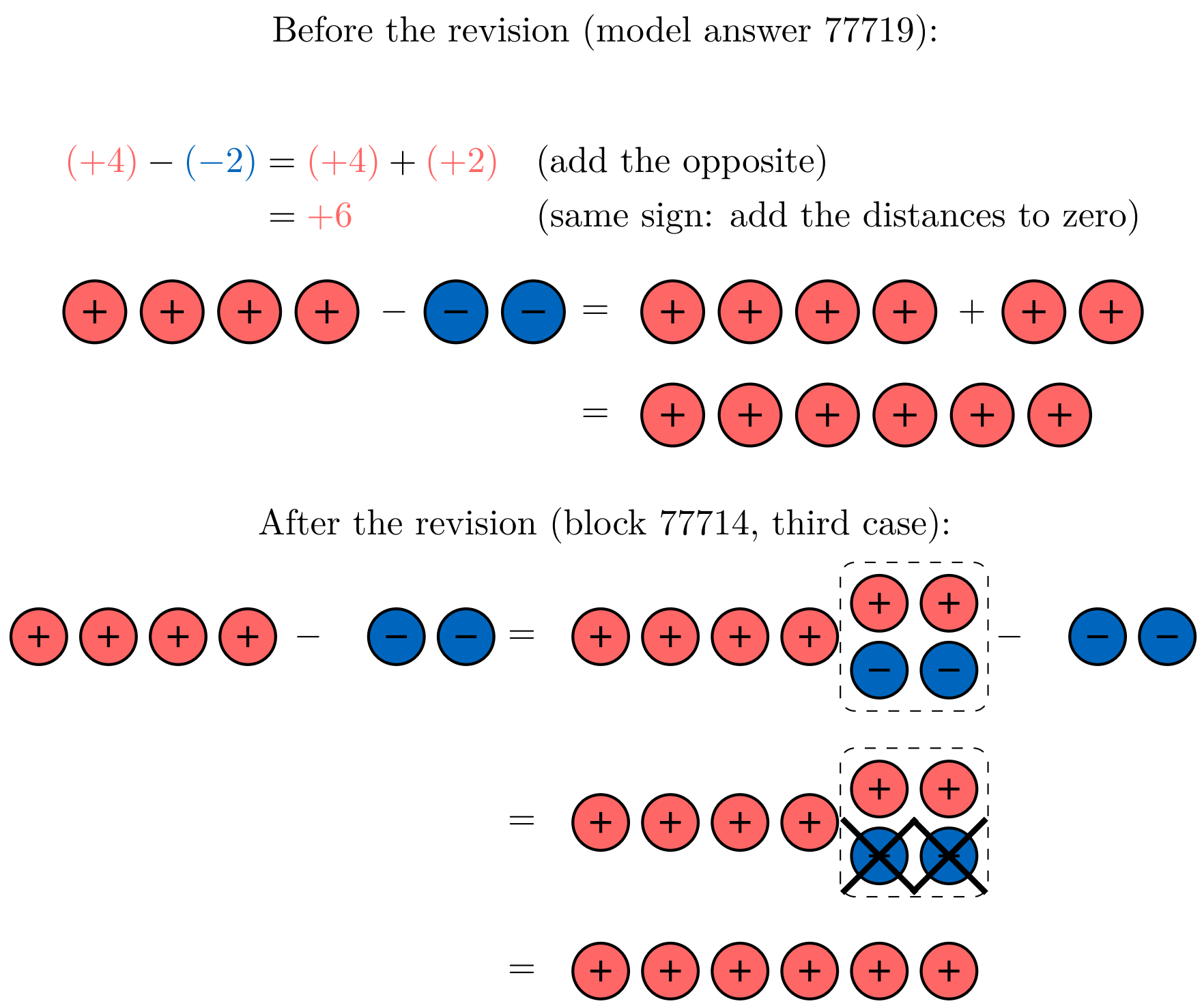

The opening problem of subtraction (Figure~\ref{fig:takeaway}) compares, with drawings, $(+3)-(+2)$ (two positives taken away from three, crossed out) with $(+3)+(-2)$ (two zero pairs cancelled), and $(-3)-(-2)$ with $(-3)+(+2)$: in each case the two operations leave the same units. In the earlier state the problem generalised from these two cases, in which the units to take away are present, to the rule ``subtracting a number is adding its opposite''. The case in which the representation is decisive was not drawn: in $(+4)-(-2)$ there is no negative unit to take away. The model answer drew $(+4)-(-2)$ directly as $(+4)+(+2)$, that is, with the rule it was supposed to justify. The revised problem adds this case (Figure~\ref{fig:subtraction}): two zero pairs are inserted, the number is still $+4$, the two negative units can now be taken away, and six positive units are left. The chain written under the drawing is the one of the 2026 examples, with integers:
\[
(+4)-(-2)=(+4)+(+2)+(-2)-(-2)=(+4)+(+2)=+6 .
\]
Taking away $(-2)$ gave the same result as adding $(+2)$. The drawing distinguishes the two gestures: inserted pairs are framed, units taken away are crossed out, and a pair that cancels is a burst, as everywhere in the lessons; on paper, where a burst cannot be drawn, the pupil circles it. Without that distinction, the same mark, a cross, would have served for removing a unit and for cancelling a pair.

\subsection{A task and its didactic variables}\label{sec:task}
The written exercise added at Grade~6 asks pupils to calculate $(+4)-(-2)$ twice, first inserting two zero pairs, then three; to explain why the starting number stays $+4$ and why both calculations give the same result; to check by adding $(-2)$ to the result; and, without drawing, to predict the result with four pairs and justify it. Its didactic variables are the signs of the two numbers, whether the units to take away are present, and the number of pairs inserted. Four strategies can be anticipated. (i)~Taking away from what is there, which gives $+2$ (two units removed from four, the sign ignored) or an impasse. (ii)~Inserting the pairs, removing, and reading $+6$, the strategy taught. (iii)~Applying ``minus minus makes plus'' without a drawing, which gives $+6$ without an explanation, and which the request to explain does not reward. (iv)~A jump on a number line, which some pupils may bring from elsewhere. The feedback is internal to the drawing: two different numbers of inserted pairs must give the same result, and the inverse addition must return the starting number. The prediction with four pairs, made without drawing, is the part of the task that the procedure alone cannot answer: it requires the statement that pairs left over leave nothing. A correct prediction alone does not establish an understanding of the invariance.

This is exactly the situation in which the literature locates the weakness of chip models. \citet{Vig2014} call it a breaking point and warn that inserting zero pairs can become a graphical procedure; in \citet{NurnbergerHaag2015PMENA}, the pooled set of problems that require inserting pairs before removing is where the chip group did worse, and the gap had not narrowed five weeks later. That study already assessed pupils' explanations of invariance under different sums of opposites; it did not assess transfer to a distributive derivation. The design does not dispute these results. It bets on what they did not measure: that the gesture of inserting a zero, made explicit and asked for in writing, is the gesture that the proof of Principle~3 reuses. If this bet is wrong, the progression pays the cost observed by Nurnberger-Haag without the benefit claimed here (Hypotheses~\ref{hyp:gesture} and~\ref{hyp:subtraction}).

\subsection{When the drawing stops}
The decimal signed numbers are introduced at Grade~6 by a definition that keeps everything: ``everything built for signed integers applies to them unchanged: opposite, distance to zero, comparison, addition and subtraction''. The revision adds one sentence: only the drawing changes; the units are no longer drawn, and the rules found with them, stated with the distance to zero, are kept. Units of a tenth could have been drawn; the choice not to draw them is made so that the rules, rather than the drawing, carry the extension.

\section{Principle 3: the product of two negatives, derived from the zero pair}\label{sec:p3}

\begin{principle}
The sign of a product is established in two steps: a positive whole number times a signed number is a repeated addition of teams, and the order of the factors does not change a product; the product of two negatives, which is not a repetition, is derived from the zero pair by keeping the distributive law. The derivation is presented as what it is: the only sign compatible with the laws of calculation, not a fact read on a model.
\end{principle}

\subsection{Repetition, and where it stops}
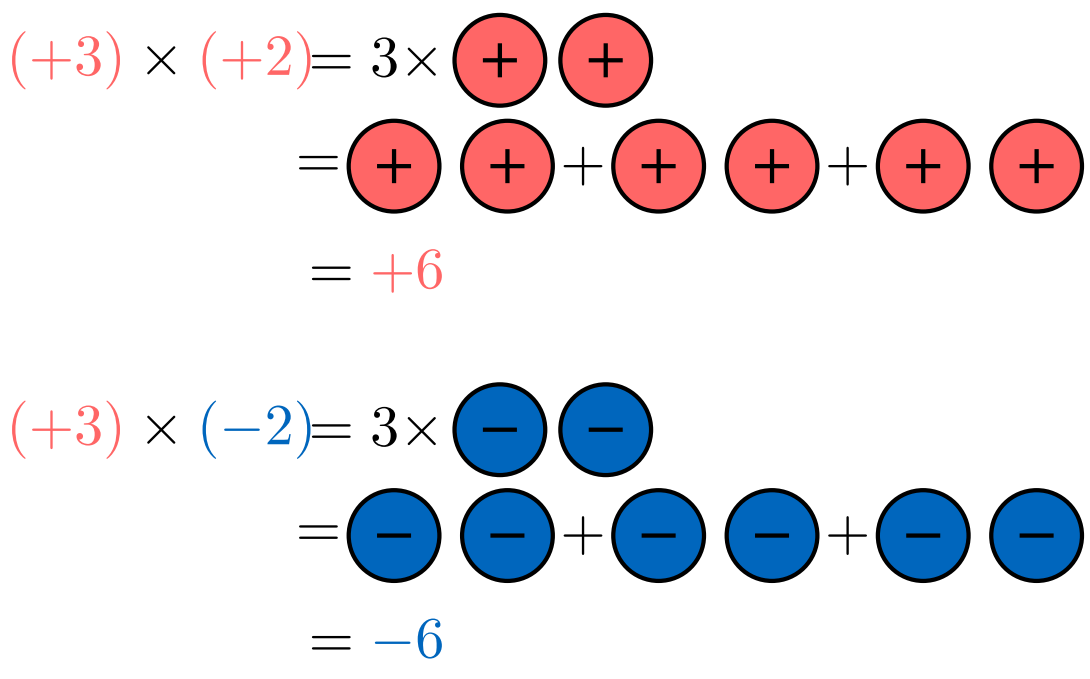
\begin{figure}[htbp]
\centering
\begin{tikzpicture}[scale=1]
\node[right] at (-0.2,0.7) {$\textcolor{colordef}{(+3)}\times \textcolor{colordef}{(+2)} $};
\node[right] at (1.8,0.7) {$=3\times $};
\tribuPositive{3.2}{0.7};
\tribuPositive{3.9}{0.7};
\node at (2,0) {$=$};
\tribuPositive{2.5}{0};
\tribuPositive{3.25}{0};
\node at (3.75,0) {$+$};
\tribuPositive{4.25}{0};
\tribuPositive{5}{0};
\node at (5.5,0) {$+$};
\tribuPositive{6}{0};
\tribuPositive{6.75}{0};
\node[right] at (1.8,-0.7) {$=\textcolor{colordef}{+6}$};
    \end{tikzpicture}

\bigskip
\begin{tikzpicture}[scale=1]
\node[right] at (-0.2,0.7) {$\textcolor{colordef}{(+3)}\times \textcolor{colorprop}{(-2)} $};
\node[right] at (1.8,0.7) {$=3\times $};
\tribuNegative{3.2}{0.7};
\tribuNegative{3.9}{0.7};
\node at (2,0) {$=$};
\tribuNegative{2.5}{0};
\tribuNegative{3.25}{0};
\node at (3.75,0) {$+$};
\tribuNegative{4.25}{0};
\tribuNegative{5}{0};
\node at (5.5,0) {$+$};
\tribuNegative{6}{0};
\tribuNegative{6.75}{0};
\node[right] at (1.8,-0.7) {$=\textcolor{colorprop}{-6}$};
    \end{tikzpicture}
\caption{Multiplication as repeated addition of teams (block 352, Grade~7): $(+3)\times(+2)$ and $(+3)\times(-2)$ as three teams of two units.}
\label{fig:teams}
\end{figure}

The Grade~7 opening problem on multiplication draws $(+3)\times(+2)$ and $(+3)\times(-2)$ as three teams of two units (Figure~\ref{fig:teams}), positive or negative, and writes $3\times(-2)=(-2)+(-2)+(-2)=-6$; $(-2)\times(+3)$ follows from the commutativity of the product. A product such as $(-3)\times(-2)$ cannot be read the same way: ``$-3$ times'' is not a number of repetitions. This is the point at which, according to Glaeser's sixth obstacle, a model built for addition is pushed into the multiplicative field where it does not apply \citep{Glaeser1981}. The units are an additive model too, and the design does not ask them to show the product.

\subsection{The derivation and what it assumes}
The same opening problem then writes, from the zero pair $(+3)+(-3)=0$:
\[
\begin{aligned}
\bigl((+3)+(-3)\bigr)\times(-2)&=0 &&\text{(the factor $(+3)+(-3)$ is zero)}\\
(+3)\times(-2)+(-3)\times(-2)&=0 &&\text{(distributivity)}\\
-6+(-3)\times(-2)&=0 &&\text{($(+3)\times(-2)=-6$)}\\
(-3)\times(-2)&=+6 &&\text{(added to $-6$ it gives $0$: it is the opposite of $-6$)}
\end{aligned}
\]
The argument uses multiplication by a whole number as repeated addition, including the empty sum ($0\times(-2)=0$: adding $(-2)$ zero times gives nothing), the uniqueness of the opposite (a number which added to $-6$ gives $0$ is the opposite of $-6$), and the extension of distributivity to signed numbers, which forces the missing product; the commutativity of the product is not needed here. In the earlier state, the problem asserted the first (``we know that $0\times(-2)=0$''), applied the second without saying it was a choice, and passed to the last line without a reason. The revised problem states that ``we want the rules of calculation known for natural numbers to stay true with signed numbers, in particular the distributive law'', justifies $0\times(-2)=0$ (adding $(-2)$ zero times gives nothing), annotates the last line, and concludes that a positive product of two negatives ``is the only sign that keeps the distributive law true''. This is Hankel's position \citep{Hankel1867} and that of the French documents \citep{Eduscol2016Relatifs,MEN2026Exemples}: the rule is not discovered in the objects, it is forced by the decision to keep the laws. Figure~\ref{fig:zero} sets the three uses of the inserted zero side by side: both arguments use a sum of opposites that represents zero, added to the number subtracted from in the subtraction, used as a factor in the derivation. This structure is not exclusive to the units; the units make it drawable.

\begin{figure}[htbp]
\centering
\[
\begin{aligned}
(+4)-(-2)&=(+4)+\underbrace{(+2)+(-2)}_{\text{a zero}}-(-2)=+6 &&\text{the subtraction drawn in Section~\ref{sec:subtraction}}\\[4pt]
3.1-(-2)&=3.1+\underbrace{2+(-2)}_{\text{a zero}}-(-2)=5.1 &&\text{the 2026 French examples}\\[4pt]
\underbrace{\bigl((+3)+(-3)\bigr)}_{\text{a zero}}\times(-2)&=0, \text{ so } (-3)\times(-2)=+6 &&\text{the product of two negatives}
\end{aligned}
\]
\caption{The same zero, three times. A zero is written in the form $a+(-a)$ without changing anything; in the first line pupils draw it as inserted zero pairs, in the third it is the factor that distributivity spreads over $(-2)$.}
\label{fig:zero}
\end{figure}

\subsection{A pattern, a conjecture, and two pathways}
A second opening problem, written for the International Baccalaureate chapters, lets the second factor decrease: $3\times2=6$, $3\times1=3$, $3\times0=0$, and continuing the pattern gives $3\times(-1)=-3$; then $2\times(-2)=-4$, $1\times(-2)=-2$, $0\times(-2)=0$, and the pattern gives $(-1)\times(-2)=+2$. The problem ends with ``a pattern suggests a rule; it does not prove it''. The two curricula use these problems differently. The French \emph{quatri\`eme} chapter includes the derivation by distributivity but not the pattern task. The IB MYP~1 chapter includes the pattern task; neither MYP~1 nor MYP~2 includes the derivation, whose intended postponement to MYP~3 is documented in design notes, not established by the audited chapters. The IB pathway thus gives pupils, for two years, a conjecture supported by regularity, which is consistent with the publisher's sequence cited in Section~\ref{sec:background}; the French pathway gives the derivation from the first year of multiplication. The revisions of the derivation (Section~\ref{sec:revisions}) therefore reach the French chapter only.

\subsection{Why not model the product}
\citet{Battista1983} presented the charge model as complete for the four operations, and a product such as $(-4)\times(-3)$ can be drawn as removing four groups of three negative units from a zero made of enough pairs. The design does not use this drawing, for two reasons. The first is that removing groups ``$-4$ times'' requires a convention that is not in the meaning of the units; the second is empirical, and is the one \citet{NurnbergerHaag2015PMENA} reports: that kind of product is among the problems on which the chip group did worse. The choice not to model the product is a design decision, not a mathematical necessity. Division is then read from multiplication, its inverse operation: since $(-3)\times(-2)=+6$, $(+6)\div(-3)=-2$; the definition states the same sign rules for a quotient, the divisor being non-zero.

\section{Principle 4: each representation for what it carries}\label{sec:p4}

\begin{principle}
Situations are kept for the encounter with negative numbers, the number line for position, order and distance, and the units for equivalence and compensation. The units are presented first and taken up next to each new notion; the hinge between units and line is the distance to zero. The drawing is abandoned when the rules can stand alone.
\end{principle}

The units do not show that $-5<-2$; the line does, by the rule ``the lower, or the further left, the smaller''. The Grade~6 method for comparing two numbers is stated with the distance to zero: of two negatives, the one closer to zero is greater. The revised definition of the distance to zero gives the hinge: for an integer, it is the number of units left once all zero pairs are cancelled, so that $+3$ and $-3$, drawn with three units each, are both three steps from zero; the revised Grade~5 method places $-4$, drawn with four negative units, four steps below zero. The order relation is known to be a productive resource for pupils \citep{Bishop2014} and has been studied in lessons dominated by the line \citep{Schindler2017}; keeping it on the line is consistent with the 2026 examples, which do the same.

The drawing is abandoned three times, each time deliberately: for decimal signed numbers, where the rules found with the units are kept (Section~\ref{sec:p2}); for the product of two negatives, which is derived, not drawn (Section~\ref{sec:p3}); and for chains of operations and expressions at Grade~7, where calculations are written without units. The Eduscol warning against relying on a single model \citep{Eduscol2016Relatifs} is not contradicted: one representation is presented first and taken up as a common thread, but it is not asked to carry what it cannot.

\section{Principle 5: the drawing first, then the chain of equalities}\label{sec:p5}

\begin{principle}
A calculation with units is first drawn as an equality between drawings, then written as a chain of equalities whose lines transcribe the drawing; the sign of a number is written in brackets and in the colour of its units, to distinguish it from the sign of the operation.
\end{principle}

The drawing macros of the chapters write each addition and subtraction as an equality between drawings, line for line, which plays for signed numbers the role that equations between bars play for fractions in the companion paper. The symbolic chain is written underneath, with the numbers in brackets and coloured like their units: $(+4)$ in red, $(-2)$ in blue. Writing the sign on the object, then in brackets, is a way of negotiating the roles of the minus sign that \citet{Vlassis2004,Vlassis2008} describe; the written chain is where these roles meet.

In the earlier state, units appeared in 50 exercise blocks of Grade~6, but in 45 answers and 5 statements: the corrections drew, and no written task required pupils to draw. Of the 33 written exercises of Grade~6 and the 38 of Grade~7, none asked for a drawing of units; the written exercises of Grade~5 asked for number lines. The revisions add twelve written exercises in which the pupil draws: four at Grade~5 (draw $-3$ and $+2$; draw $-2$ in two ways and circle the zero pairs; draw zero without units and with units; circle the zero pairs of a drawing and write the number left) and eight at Grade~6 (four additions and four subtractions, each drawn, then written as a chain). Their answers are drawn as the lessons draw, with bursts for the pairs that cancel; the first exercise that uses pairs tells the pupil that a burst is circled on paper. They are graded on a photograph or on a digital slate; their criteria accept any layout and a drawing in one colour if every sign is legible, distinguish the circled pair from the crossed unit, and penalise units taken away that were never drawn and pairs that change the starting number.

The annotation of chains (Table~\ref{tab:corpus}) shows the state of the exercise answers: at Grade~6, 28 of the 68 multi-line chains carried a side note before the revisions, of which 21 a verbal reason (``add the opposite'', ``sum of signed numbers''); the others noted a side calculation such as $(+3)+(-2)=+1$. The eight new answers add five chains, all with verbal reasons. At Grade~7, 48 of 94 chains carry a note. In the lessons, every multi-line chain of Grade~6 is annotated; at Grade~7, the revisions annotate the three chains of the opening problem on multiplication that had no note (``optional $+$ sign'', ``repeated addition'', ``commutativity'') and the last line of its derivation.

\section{The progression from Grade~5 to Grade~7}\label{sec:progression}

\paragraph{Grade 5: the units, then the situations, then the line.} The lesson opens with the particles and the definition of the integers (Section~\ref{sec:p1}), then meets negative numbers in three situations, reads them with a method that adds the row of units, places them on a line through zero, reads heights on a cliff and quiz scores on a line, compares and orders values in context, and ends with two remarks: the numbers met are integers, and opposite numbers on the line are at the same distance from zero. Its exercises count the units of a drawing, write integers from words, ask the pupil to draw integers (four written exercises), write situations as numbers, read lines through zero, order values in context (written, on a drawn line), and, after the line, find opposites, with answers that begin with the units and end on the line.

\paragraph{Grade 6: adding and subtracting with units, the line, the decimals.} The section ``What are integers?'' links the Grade~5 blocks and adds a remark on the families of numbers (natural numbers, integers, and numbers such as $2.5$ that are neither). The number line and the distance to zero, comparing and ordering, addition, subtraction (with the decisive case), the distance between two points read as a difference, the decimal signed numbers, and chains of additions and subtractions follow. Its exercises include the two written groups in which the pupil draws, at the head of the addition and subtraction sections.

\paragraph{Grade 7: products and quotients.} The first five sections link the Grade~6 lessons, so that a pupil who meets Grade~7 first finds the same representation. Multiplication (repetition, then the zero pair and distributivity, then the pattern), the sign of a product of several factors, powers of negative numbers, division read from multiplication, decimal signed numbers, the order of operations, the vocabulary of expressions and the calculator follow.

\paragraph{Projections.} Table~\ref{tab:projections} lists the five chapters assembled from the progression. The French \emph{cinqui\`eme} chapter shows Grade~5 and Grade~6, with the unit of the particles first and the unit of the situations second; the PYP~5 chapter shows Grade~5 after a unit on degrees Fahrenheit. That unit does not teach negative numbers, but its thermometers are graduated down to $-10$ on the Celsius scale: in PYP~5 the units are the first explicit introduction of the integers, not their first visual occurrence. The earlier IB plan had kept the particles for MYP~1, because the PYP~5 scope does not include the addition of integers, which the particles' opening problem already writes with numbers; the revision brings them to PYP~5 deliberately, for the sake of one representation from the first explicit introduction.

\begin{table}[htbp]
\centering\small
\caption{Chapters assembled from the progression (link tables of 28~October~2026). ``Links'' counts the blocks a chapter displays, lessons and exercises. The link tables changed for \emph{cinqui\`eme}, PYP~5 and MYP~1; the other chapters receive the revised blocks through their links.}
\label{tab:projections}
\begin{tabular}{llll}
\toprule
Curriculum & Chapter & Source levels & Links, before $\to$ after \\
\midrule
French, \emph{cinqui\`eme} & Signed numbers & 5, 6 & $225\to237$ \\
French, \emph{quatri\`eme} & Multiplying and dividing & 5, 6, 7 & $330$ \\
IB, PYP 5 & Temperature and negative numbers & 5, 6 & $108\to120$ \\
IB, MYP 1 & Positive and negative numbers & 5, 6, 7 & $117\to125$ \\
IB, MYP 2 & Integers and signed numbers & 6, 7 & $143$ \\
\bottomrule
\end{tabular}
\end{table}

\section{Anticipated misconceptions}\label{sec:misconceptions}

Each misconception below is documented in the literature or in the materials' own warnings; the design response is stated with the block that carries it.

\begin{enumerate}[leftmargin=*]
\item \emph{``Minus and minus make plus'' applied to addition}, $(-3)+(-2)=+5$. The drawing shows five negative units and no pair to set aside; the rule of addition is stated with the distance to zero before any rule of products.
\item \emph{Adding always makes larger, subtracting always smaller} \citep{Bishop2014}. A remark states that the sign of the number added or subtracted decides, with $(+5)+(-3)=+2$ and $(+5)-(-3)=+8$.
\item \emph{Taking away more than one has} \citep{Bishop2014}: $(+2)-(+5)$ read as impossible, or as $3$. The written exercise asks for at least three zero pairs to be inserted before five positive units are taken away.
\item \emph{Comparing negative numbers by their magnitude}, $-5>-2$. Comparison is taught on the line and with the distance to zero, not on the units (Section~\ref{sec:p4}); a collection with more negative units is not ``more''.
\item \emph{Two signs side by side}, $3+-5$. A remark states that two signs never stand side by side and shows the brackets or the simplified writing $3-5$.
\item \emph{The square of a negative number}, $-3^2$ read as $(-3)^2$. A remark on powers states that, without brackets, the exponent applies to $3$ only: $-3^2=-(3^2)=-9$, whereas $(-3)^2=9$.
\item \emph{A state read as a change}, a temperature of $-4\,^\circ$C read as a fall of four degrees \citep{Vergnaud1982}. The Grade~5 method reads a quantity below a zero, the units count the steps from that zero; changes of temperature appear at Grade~6, as differences.
\item \emph{A graphical procedure without meaning} \citep{Vig2014}: inserting zero pairs because one is told to. The written exercises ask why the number does not change and to predict without drawing (Section~\ref{sec:task}).
\end{enumerate}

\section{Discussion}\label{sec:discussion}

\subsection{Position with respect to existing recommendations}
The progression's justifications are compatible with the examples that accompany the 2026 French programme \citep{MEN2026Exemples}: the opposite through a zero sum, subtraction through an inserted zero, the product of two negatives through distributivity, the line for order. Those examples do not prescribe units, and they also use the line for the opposite; the progression gives their steps a drawable form for pupils of ten to thirteen. It shares the diagnosis of the 2016 resources and of the two French groups that concrete additive models may hinder the product \citep{Eduscol2016Relatifs,IREMAquitaine2008,Matheron2016}, but it differs from their remedies (equalities with a missing term, calculation programmes) by keeping a concrete representation for the additive part. It belongs to the second of Janvier's families, the models of opposite objects \citep{Janvier1985}, and inherits their known difficulty with subtraction \citep{Janvier1985,Vig2014,NurnbergerHaag2015PMENA}.

\subsection{What the audit revealed}
The audit of the state of 28~October~2026, before the revisions, found five discrepancies between the stated design and the documents.
\begin{enumerate}[leftmargin=*]
\item \emph{Two first meanings.} Grade~5 entered by situations and the number line and contained no unit (15 of its 18 lesson blocks cited a situation); it defined the opposite by the distance to zero. Pupils whose chapter began at Grade~5 met a first meaning (a level below zero) and, one level later, a second (a collection of units).
\item \emph{No written task required pupils to draw.} Units appeared in 50 exercise blocks of Grade~6, of which 45 answers and 5 statements; none of the 75 written exercises of Grades 5 to 7 asked for a drawing of units. The corpus says nothing of what teachers ask in class.
\item \emph{The decisive subtraction was not drawn.} The rule ``subtracting is adding the opposite'' was generalised from two cases in which the units to take away are present; the case that requires inserting zero pairs was drawn with the rule itself (Section~\ref{sec:subtraction}).
\item \emph{The derivation did not state its assumptions.} The product of two negatives was derived by distributivity, but the text asserted $0\times(-2)=0$, did not say that keeping distributivity is a choice, and passed to the result without naming the opposite (Section~\ref{sec:p3}).
\item \emph{The change of domain was not announced.} The decimals inherited ``everything'' from the integers, without saying that the units were no longer drawn (Section~\ref{sec:p2}).
\end{enumerate}
The audit also revealed errors in the author's first counts, found by the external reading (Section~\ref{sec:method}): a chain count of zero annotated answers at Grade~6, where there were 28, came from a pattern that missed one of the two written forms of a margin note. Counts made by script need the same scrutiny as the statements they support.

\subsection{Revisions made after the audit}\label{sec:revisions}
The revisions were written as drafts, read twice by the external reviewer, corrected, approved by the author, and merged into the source documents on 28~October~2026; the chapters were regenerated from the revised sources and published on the public site the same day. They are listed here because the analysis of the preceding sections describes the revised state, and because a reader of the earlier chapters will meet the earlier state.
\begin{enumerate}[leftmargin=*]
\item \emph{One representation from the first lesson.} The particles' opening problem, the definition of the integers and the video that introduces them moved from Grade~6 to the head of Grade~5, with their identifiers, and are shown at Grade~6 and Grade~7 by links; their text was retouched (each particle stands for one unit; opposites ``drawn without zero pairs''; the number lines that preceded the teaching of the line were removed). The Grade~5 definition by situations became a method, the definition of opposites on the line became a remark beginning with the units, and a remark that redefined the integers was reworded. Three groups of exercises moved down with them.
\item \emph{The units next to the line.} With a step of~1, the number of steps from zero is the number of units; the distance to zero is the number of units left once the zero pairs are cancelled; the answers on opposites begin with the units.
\item \emph{The decisive subtraction}, drawn with inserted pairs (Figure~\ref{fig:subtraction}).
\item \emph{Twelve written exercises in which the pupil draws}, four at Grade~5 and eight at Grade~6, with criteria for grading a photograph.
\item \emph{The derivation states its assumptions} (Section~\ref{sec:p3}), and its first chains are annotated.
\item \emph{The change of domain is announced} for the decimals.
\end{enumerate}
In the chapters, the French \emph{cinqui\`eme} now begins with the particles and meets the situations second, and the PYP~5 chapter opens its part on negative numbers with them (Table~\ref{tab:projections}). The video that opened the Grade~5 situations, already published, shows situations only; it now opens the second unit, and the video of the particles opens the first.

\subsection{What the progression does not address}
Grade~7 introduces negative fractions through notation and quotients ($-\tfrac12$; $1\div(-3)=-\tfrac13$); it does not provide a systematic progression for the arithmetic of rational numbers, which belongs to the last stage of the fraction progression analysed in the companion paper. It does not state that multiplying by $-1$ gives the opposite, which the 2026 examples do. The meaning of the minus sign as an operator on a letter ($-x$), central in algebra \citep{Vlassis2004}, and the use of signed numbers as displacements or vectors are outside its scope. The order of the integers is taught on the line, without a definition on the units.

\subsection{Hypotheses for an empirical evaluation}
The design rests on hypotheses that the documents cannot test.
\begin{hypothesis}[The gesture]\label{hyp:gesture}
Pupils who have inserted zero pairs in subtraction explain the derivation of $(-3)\times(-2)=+6$ (why the zero pair is used, why $-6+x=0$ gives $x=+6$) better than pupils taught the same derivation without the unit representation.
\end{hypothesis}
\begin{hypothesis}[Subtraction]\label{hyp:subtraction}
On subtractions that require inserting zero pairs, pupils taught with the units do worse than pupils taught with a path model, as in \citet{NurnbergerHaag2015PMENA}, where the gap did not narrow but widened at five weeks; the hypothesis is that written explanation and prediction tasks, absent from that study, prevent this widening. It is refuted if the delayed gap is not smaller than the immediate one.
\end{hypothesis}
\begin{hypothesis}[One representation]
Presenting the units before the situations at Grade~5 does not impair reading and placing negative numbers in situations, and pupils coordinate a state (a temperature) with its units without reading it as a change.
\end{hypothesis}
\begin{hypothesis}[Order]
Pupils taught with the units compare negative numbers on the line without judging a collection with more negative units to be larger.
\end{hypothesis}
\begin{hypothesis}[Written drawings]
Drawings of units photographed by pupils can be graded with the published criteria with an agreement between graders comparable to that of numerical answers.
\end{hypothesis}
Hypothesis~\ref{hyp:subtraction} is the one on which the existing evidence weighs against the design; Hypothesis~\ref{hyp:gesture} is the one on which the design depends. A study that tested only the first would measure the cost of the choice without its intended benefit. A test of both would compare three groups with the same duration and the same teaching of the derivation (ordinary chips, units with the written explanation tasks, and a path model), with a pre-test of opposite sums and new items graded blind. Outcomes would distinguish arithmetic accuracy, explanations of invariance and transfer to unfamiliar distributive arguments, assessed immediately and five weeks later; non-inferiority margins and inter-rater agreement criteria would be specified in advance.

\subsection{Limitations}
No data on pupils are reported. The author designed the progression, made the audit and wrote the revisions; the external reading was done by a language model, not by human experts in the didactics of signed numbers. Counts concern the English halves of the source documents; the list of situation words and the criteria for chains are explicit but conventional. The link tables establish the content and order of the chapters, not what pupils experienced. Several references are known only at the level of their abstract or through secondary sources, and this is stated where it matters: Glaeser's list of obstacles is taken from three concordant secondary sources, and the direction of Liebeck's comparison with a number-line group was not available to the author.

\section{Conclusion}\label{sec:conclusion}

The progression analysed here presents one representation of signed numbers first, positive and negative units that cancel in pairs, and takes it up at each notion until the rules can stand without it. Its central claim is narrower than ``units give more meaning than the line'': the units make visible that an integer has many drawings, and they give the pupil the gesture of inserting a zero, which subtraction and the derivation of the sign rule share. The product of two negatives is not modelled; it is derived, and the derivation says that it holds because the laws of calculation are kept. The audit showed that the documents did not realise this design at its first level, in its exercises and in its decisive case, and the revisions made the representation the first one met, asked pupils to draw, and drew the case on which the representation's value depends. Whether that case is the representation's strength, as the design hopes, or its weakness, as a classroom comparison suggests, is an empirical question that this paper can only state.

\section*{Competing interests}
The author is the founder of CommeUnJeu, which publishes the materials analysed in this paper and offers paid services on its website.

\section*{Acknowledgements}
The author used AI language tools while preparing this paper: Claude (Anthropic) assisted with the audit script, the drafting of the revisions and of the paper, and the search for and checking of bibliographic references, and ChatGPT (OpenAI) assisted with the critical reading of the plan and of the revisions. The author designed the progression and takes full responsibility for the content.

\bibliographystyle{abbrvnat}
\bibliography{relatifs_refs}

\appendix
\section{Public resources}\label{app:resources}

The chapters are served at the addresses below, in English and in French. All addresses begin with \texttt{https://commeunjeu.org/cours/}; in the French slugs, \texttt{grade-}$N$ denotes the $N$-th year of French compulsory schooling after kindergarten (CP is 1), not a level of the progression. They show the revisions of 28~October~2026, published the same day (pages checked on 28~October~2026).
\begin{center}\small
\begin{tabular}{ll}
\toprule
French, \emph{cinqui\`eme} & \texttt{fr-grade-7-signed-numbers/} \\
French, \emph{quatri\`eme} & \texttt{fr-grade-8-signed-numbers-multiplication-and-division/} \\
IB, PYP 5 & \texttt{ib-pyp-5-temperature-and-negative-numbers/} \\
IB, MYP 1 & \texttt{ib-myp-1-positive-and-negative-numbers/} \\
IB, MYP 2 & \texttt{ib-myp-2-integers-and-signed-numbers/} \\
\bottomrule
\end{tabular}
\end{center}
Each chapter page embeds the videos attached to the sections it uses, published on two channels, one per language: \texttt{youtube.com/@CommeUnJeuEN} (English) and \texttt{youtube.com/@CommeUnJeu} (French).

\section{Audit conventions}\label{app:audit}

\begin{enumerate}[leftmargin=*]
\item \textbf{Corpus.} The six source documents in their two dated states of 28~October~2026, before and after the revisions. Counts use the English half of each block; a link displays a block written elsewhere and is counted where it is written, and separately as a link.
\item \textbf{Units.} A block draws units when its drawing uses the macros of the positive and negative units, of the drawn additions and subtractions, or of the cancelled pair.
\item \textbf{Situations.} A block cites a situation when its English text contains, as a whole word, an item of a closed list of about twenty (temperature, thermometer, degrees Celsius, sea level, height, depth, diver, floor, lift, bank, debt, account, altitude, score, quiz, colder, warmer and a few variants); the list selects 15 of the 18 lesson blocks of Grade~5 before the revisions (all but the video, the remark on integers and the definition of opposites).
\item \textbf{Chains.} A chain is a multi-line calculation (\texttt{aligned} or \texttt{align}, at least two lines). It carries a side annotation when a line has anything after the alignment marker of its margin; a verbal reason when that note contains words; it is fully annotated when every line but the last does. These are counts of annotation practice, not judgements that each step is justified.
\item \textbf{Exercises.} Exercises and multiple-choice questions are counted apart; a written exercise is one graded from a photograph or a digital slate.
\end{enumerate}
An independently inspectable reanalysis would need the dated source files, the script and its item-level output; they are held by the author and are not provided with this manuscript.

\end{document}